\documentclass[11pt]{amsart}

\usepackage{amsxtra}
\usepackage{amsfonts}
\usepackage[latin1]{inputenc}
\usepackage{epsfig}
\usepackage[T1]{fontenc}
\usepackage{amsmath,amssymb,latexsym,amsthm}
\usepackage[all]{xy}
\usepackage{mathrsfs}
\usepackage{psfrag}

\newtheorem{theorem}{Theorem}[section]
\newtheorem{lemma}[theorem]{Lemma}

\newtheorem{definition}[theorem]{Definition}
\newtheorem{proposition}[theorem]{Proposition}

\newtheorem{remark}[theorem]{Remark}
\newtheorem{example}[theorem]{Example}
\newtheorem{notation}[theorem]{Notation}
\newtheorem{observation}[theorem]{Observation}

\numberwithin{equation}{section}

\CompileMatrices

\begin{document}

\title{The cyclic bar construction on $A_\infty$ $H$-spaces}

\author{Vigleik Angeltveit} \thanks{This research was partially conducted during the period the author was employed by the Clay Mathematics Institute as a Liftoff Fellow}

\maketitle

\newcommand{\A}{\mathcal{A}}
\newcommand{\B}{\mathcal{B}}
\newcommand{\C}{\mathcal{C}}
\newcommand{\D}{\mathcal{D}}
\newcommand{\E}{\mathcal{E}}
\newcommand{\F}{\mathbb{F}}
\newcommand{\I}{\mathcal{I}}
\newcommand{\K}{\mathcal{K}}
\newcommand{\M}{\mathcal{M}}
\newcommand{\R}{\mathbb{R}}
\newcommand{\V}{\mathcal{V}}
\newcommand{\W}{\mathbb{W}}
\newcommand{\Z}{\mathbb{Z}}
\newcommand{\sma}{\wedge}
\newcommand{\lar}{\longrightarrow}
\newcommand{\sar}{\rightarrow}
\newcommand{\colim}{\varinjlim}
\newcommand{\hEn}{\widehat{E(n)}}
\newcommand{\bD}{\mathbf{\Delta}}
\newcommand{\noco}{\mathbf{\Delta} \Sigma}
\newcommand{\Top}{\mathcal{T}\!\mathit{op}}


\def\arrow#1{\overset{#1}{\lar}}

\begin{abstract}
We set up a general framework for enriching a subcategory of the category of noncommutative sets over a category $\C$ using products of the objects of a non-$\Sigma$ operad $P$ in $\C$. By viewing the simplicial category as a subcategory of the category of noncommutative sets in two different ways, we obtain two generalizations of simplicial objects. For the operad given by the Stasheff associahedra we obtain a model for the $2$-sided bar construction in the first case and the cyclic bar and cobar construction in the second case. Using either the associahedra or the cyclohedra in place of the geometric simplices we can define the geometric realization of these objects.
\end{abstract}

\section{Introduction}
In \cite{St63} Stasheff introduced the notion of an $A_\infty$ $H$-space in terms of a family of polyhedra $\{K_n\}_{n \geq 0}$ which have become known as the Stasheff associahedra. He also defined the bar construction $BA$ for an $A_\infty$ $H$-space $A$ as a quotient of $\coprod K_{n+2} \times A^n$, and showed that this construction has the same good properties as the bar construction on a strictly associative $H$-space.

Stasheff's construction is somewhat ad hoc, and while it is not hard to generalize his construction to define the $2$-sided bar construction $B(M,A,N)$, defining the cyclic bar construction $B^{cy}(A;M)$ for an $A_\infty$ $H$-space $A$ and an $A$-bimodule $M$ requires a new ingredient.

The new ingredient is another family of polyhedra $\{W_n\}_{n \geq 0}$, introduced by Bott and Taubes in \cite{BoTo94} and named cyclohedra by Stasheff in \cite{St95}. Given an $A_\infty$ $H$-space $A$ and an $A$-bimodule $M$ we can then define $B^{cy}(A;M)$ as a quotient of $\coprod W_{n+1} \times M \times A^n$. We can also define what we call the cyclic cobar construction $C_{cy}(A;M)$ as a subspace of $\prod Hom(W_{n+1} \times A^n,M)$. These constructions can be thought of as the space-level versions of Hochschild homology and cohomology.

Given an $A_\infty$ $H$-space $A$, we can think of the bar construction $BA$ as the geometric realization of a ``simplicial'' space. Here we put simplicial in quotes because the simplicial identities hold only up to homotopy (and higher homotopies) and as a result we also need a more general notion of geometric realization.

In this paper we make the notion of a simplicial space up to homotopy precise by defining a new category enriched over spaces with the same objects as the simplicial indexing category but more morphisms. In fact, using two different descriptions of the simplicial indexing category we find two ways to do this. The first construction makes the $2$-sided simplicial bar construction $B_\bullet(M,A,N)$ for an $A_\infty$ $H$-space $A$ into a generalized simplicial space, while the second construction also makes the cyclic bar construction into a generalized simplicial space.

While it is possible to define the cyclic bar construction on an $A_\infty$ $H$-space in a similarly ad hoc way as Stasheff's original definition of $BA$, the constructions become much more transparent when put into a more general framework. Developing such a framework is the main purpose of this paper.


We will use the constructions in this paper in \cite{AnTHH} to calculate topological Hochschild homology and cohomology of certain $A_\infty$ ring spectra. Defining $THH$ in terms of the associahedra and cyclohedra helps clarify how $THH$ depends on the $A_\infty$ structure. To retain homotopical control over these generalized (co)simplicial objects, we were led to the definition of an enriched Reedy category and the appropriate generalization of the Reedy model category structure on functors from an enriched Reedy category to an enriched model category. We will return to this in \cite{An_R}.

\subsection*{Organization} In section \ref{operads} we recall the definition of an operad in a symmetric monoidal category $\C$ and some objects and structures related to operads and non-$\Sigma$ operads, including right modules. The reader who is familiar with operads might want to skip this section, referring back to it as necessary.

In section \ref{subnoco} we consider a subcategory $\A$ of the category $\noco$ of noncommutative sets, and given a non-$\Sigma$ operad $P$ in $\C$ we define a new category $\A_P$ enriched over $\C$ using the objects in $P$. The most productive way of thinking about the category $\noco$ from our point of view is as having finite sets as objects, and a morphism is a map of sets $f:S \sar T$ together with a linear ordering of each inverse image $f^{-1}(t)$. This way it is clear that replacing a point $f:S \sar T$ in $Hom_{\noco}(S,T)$ with $\bigotimes_{t \in T} P(f^{-1}(t))$ makes sense.

We are interested in two different embeddings of $\bD^{op}$ into $\noco$. The first is given by identifying $\bD^{op}$ with a subcategory of $\bD$ we call $^{01}\!\bD$ where the objects are sets of cardinality at least $2$ and the maps are order-preserving maps which preserve the minimal and maximal element. The second is given by looking at a certain subcategory of Connes' cyclic category $\bD C$, namely the category of based cyclically ordered sets and basepoint-preserving maps, which we will denote by $^0\!\bD C$.

Thus we get two generalizations of simplicial objects, either as a functor from $^{01}\!\bD_P$ to some $\C$-category $\D$, or as a functor from $^0\!\bD C_P$ to $\D$. If $P$ is the associative operad we get simplicial objects in both cases, but otherwise these two constructions are different. If $A$ is a $P$-algebra and $M$ and $N$ are right and left $A$-modules, respectively, then we can generalize the $2$-sided bar construction and make a functor $^{01}\!\bD_P \sar \D$. The second generalization of the simplicial category gives a generalization of the cyclic bar construction. If $A$ is a $P$-algebra and $M$ is an $A$-bimodule we can make a functor $^0\!\bD C_P \sar \D$ sending a cyclically ordered set $S$ with basepoint $0$ to $M \otimes A^{\otimes S-\{0\}}$ which generalizes the cyclic bar construction. Dually, we can make a functor $(^0 \! \bD C_P)^{op} \sar \D$ sending $S$ to $Hom(A^{\otimes S-\{0\}},M)$.

If we have functors $F: \A_P \sar \D$ and $R : \A_P^{op} \sar \C$ we can form the coequalizer $R \otimes_{\A_P} F$. If $P$ is the associative operad and $R$ sends an $n$-element set to $\Delta^{n-1}$ this gives geometric realization in the usual sense, and with appropriate choices for $P$ and $R$ this still gives a good model for geometric realization.

We then devote section \ref{asscyclo} to studying the Stasheff associahedra and the cyclohedra in some detail. We denote the associahedra operad by $\K$ and the cyclohedra by $\mathcal{W}$, and observe that the cyclohedra give a functor $\mathcal{W} : \bD C_\K^{op} \sar \Top$.

Finally, in section \ref{tracecotrace} we identify the kind of structure the cyclic bar or cobar construction corepresents or represents. Given a right $P$-module $R$, Markl \cite[Definition 2.6]{Markl} defines an $R$-trace on a $P$-algebra $A$ in $\D$ with target some other object $B$ in $\D$. We generalize this by letting $R$ be a functor ${}^0 \! \bD C_P^{op} \sar \C$, and a trace is now defined on a pair $(A,M)$ consisting of a $P$-algebra $A$ and an $A$-bimodule $M$. An $R$-trace is corepresented by the cyclic bar construction while an $R$-cotrace is represented by the cyclic cobar construction.

\subsection*{Acknowledgements}
An earlier version of this paper formed parts of the author's PhD thesis at the Massachusetts Institute of Technology under the supervision of Haynes Miller.

\section{Operads} \label{operads}
In this section we recall some things about operads that we will need later. We focus mostly on non-$\Sigma$ operads, which we will simply call operads. We will use the term $\Sigma$-operad for an operad in the usual sense, the few places where we will need operads with symmetric group actions. The original reference for operads is \cite{Ma72}, see also \cite{BoVo73}; for a more modern introduction the reader can see for example \cite{McSm04}. See also \cite{MSS} for a comprehensive treatment of many topics related to operads. Our formalism is inspired by \cite{Ching05}, though we focus mostly on non-$\Sigma$ operads.

One key difference from the approach in \cite{MSS} and \cite{Ching05} is that our operads have a zeroth space. We need a zeroth space in our operads to make sense of using an operad in $\C$ to enrich some category of sets over $\C$. Indeed, when the category of sets is equivalent to $\bD^{op}$, the maps coming from composition with a nullary operation should be thought of as degeneracy maps. One disadvantage is that, at least with the standard model category structure on operads in spaces, cofibrant operads are very big, and the Stasheff associahedra operad is not cofibrant.

\subsection*{Sequences and symmetric sequences}
Let $\C$ be a closed symmetric monoidal category with all countable products and coproducts. For the basic definitions $\C$ does not have to be closed, though it is a convenient technical assumption, see Remark \ref{quasi_monoidal}. Our main examples are the categories of spaces and based spaces, which we denote by $\Top$ and $\Top_*$. By a space we mean a compactly generated weak Hausdorff space, and we will assume that all our based spaces are well pointed, i.e., that $* \sar X$ is a cofibration.  Throughout the paper it makes sense to use the category of (based) simplicial sets instead. We will denote the monoidal structure by $\otimes$ and the unit object by $I$. We will assume that $\varnothing \otimes A=\varnothing$ for all $A$, where $\varnothing$ is the initial object. This is automatic if $\C$ is pointed (\cite[Proposition 1.13]{Ching05}). If $\C$ is a specific category, we will revert to the usual notation for that category.

Let $\D$ be a symmetric monoidal $\C$-category \cite[Definition 1.10]{Ching05}, also with all countable products and coproducts. This means that $\D$ is enriched, tensored and cotensored over $\C$ and that the usual coherence relations between these structures are satisfied. Some constructions only require that $\D$ is tensored over $\C$, or enriched and tensored, or enriched and cotensored. The requirement that $\C$ and $\D$ have all countable limits and colimits can also be relaxed for some of our constructions.

We let $\bD$ denote the category of finite, nonempty, totally ordered sets and order preserving maps, and we let $\bD_+$ be the category of all finite totally ordered sets, i.e., $\bD$ together with the empty set. We define a sequence in $\C$ as a functor
\begin{equation} 
P : iso(\bD_+) \lar \C.
\end{equation}
The isomorphisms are required to be order-preserving, so there are no nontrivial automorphisms.

We want to compare operads and simplicial objects using the category $\bD$, so a warning about the notation is in order. When considering a simplicial object $X$ it is customary to denote $X(\{0,1,\ldots,n\}) \cong X(\{1,2,\ldots,n+1\})$ by $X_n$, while one usually writes $P(n)$ for $P(\{0,1,\ldots,n-1\}) \cong P(\{1,2,\ldots,n\})$. We will use both of these conventions for functors from finite totally ordered sets, together with some nonstandard conventions for other categories of sets, and we hope this will not cause confusion.

We define a symmetric sequence in $\C$ as a functor
\begin{equation}
P : iso(\mathcal{F}) \lar \C,
\end{equation}
where $\mathcal{F}$ is the category of finite sets. In this case the symmetric group $\Sigma_n$ acts on $P(n)$. We think of a (symmetric) sequence both as a functor from some category of finite sets and as a collection of objects indexed by the natural numbers. We will find it convenient to consider the category of all finite (totally ordered) sets when writing down coherence conditions, while (implicitly) choosing a skeleton category when taking limits and colimits. The main advantage of working with arbitrary finite sets is that we avoid constant relabeling when describing, for example, the coherence conditions the structure maps for an operad have to satisfy.

To ease the notation, given a sequence $P$ and a map $f:S \sar T$ in $\bD_+$ we will write $P[f]$ for $\bigotimes_{t \in T} P(f^{-1}(t))$. This notation is inspired by \cite[Definition 1.53]{MSS}.

\begin{definition} \label{comp_prod}
Given sequences $P$ and $Q$ in $\C$, their composition product, which we denote by $P \circ Q$, is given by
\begin{equation}
(P \circ Q)(S)=\coprod_{[f:S \rightarrow T]} P(T) \otimes Q[f]
\end{equation}
for a finite totally ordered set $S$, where the coproduct runs over all isomorphism classes of totally ordered sets $T$ and all order-preserving maps $S \sar T$.

The composition product on symmetric sequences is defined similarly, using unordered sets.
\end{definition}

\begin{remark} \label{quasi_monoidal}
This does not quite define a monoidal product unless the monoidal product in $\C$ distributes over coproducts, as the operation $\circ$ on sequences in $\C$ is not quite associative up to natural isomorphism. We refer the reader to \cite{Ching3}, which explains how to get around this by defining $M \circ N \circ P$, which maps to both $(M \circ N) \circ P$ and $M \circ (N \circ P)$. We will ignore this technicality in this paper, since the monoidal product does distribute over coproducts when $\C$ is closed and all the symmetric monoidal categories we consider are closed.
\end{remark}

\subsection*{Operads and modules}
Let $\I$ be the sequence with $\I(1)=I$ and $\I(n)=\emptyset$ for $n \neq 1$, and note that $P \circ \I \cong P \cong \I \circ P$.

\begin{definition}
An operad is a sequence $S \mapsto P(S)$ together with a unit map $\I \sar P$ and an associative and unital map
\begin{equation}
P \circ P \lar P.
\end{equation}
A right $P$-module is a sequence $R$ together with an associative and unital map
\begin{equation}
R \circ P \lar R,
\end{equation}
and a left $P$-module $L$ is a sequence together with an associative and unital map
\begin{equation}
P \circ L \lar L.
\end{equation}
A $\Sigma$-operad is defined similarly.
\end{definition}

\begin{remark}
Suppose $P$ is a sequence in $\C$ and $L$ is a sequence in $\D$. Then we can define the composition product $P \circ L$ as above, and this is again a sequence in $\D$. If $P$ is an operad in $\C$, the above definition of a left $P$-module still makes sense if $L$ is a sequence in $\D$ rather than in $\C$.

The same holds for right $P$-modules.
\end{remark}

Thus an operad structure on $P$ is a collection of maps $P(T) \otimes P[f] \sar P(S)$ for each $f : S \sar T$ in $\bD_+$ satisfying certain conditions.

\begin{definition} \label{circle_notation}
Let $s \in S$, and let $S \cup_s T$ denote $(S-\{s\}) \coprod T$ ordered in the obvious way, with $s_1 < t < s_2$ if $s_1<s<s_2$. Let $f:S \cup_s T \sar S$ be the map which restricts to the identity on $S-\{s\}$ and sends $T$ to $s$. Let $\circ_s : P(S) \otimes P(T) \rightarrow P(S \cup_s T)$ be the composite $P(S) \otimes P(T) \sar P(S) \otimes P[f] \sar P(S \cup_s T)$, where the first map is given by $I \sar P(1)$ on each of the $P(1)$-factors in $P[f]$.
\end{definition}

When describing the structure maps for an operad, it is enough to give the maps $\circ_s$, see the discussion about pseudo-operads in \cite[\S 1.7.1]{MSS}. 

We embed the category $\D$ in the category of sequences in $\D$ by letting $A(S)=A$ if $S=\varnothing$ and $A(S)=\varnothing$ otherwise. (It is also possible to let $A(S)=A$ for all $S$.)

We say that a $P$-algebra structure on $A$ is a left $P$-module structure on the corresponding sequence. We find that $A[f]=A^{\otimes T}$ for $f : S \sar T$ if $S=\varnothing$, and $A[f]=\varnothing$ otherwise. The maps $P(T) \otimes A[f] \sar A(S)$ reduce to maps $P(n) \otimes A^{\otimes n} \sar A$ satisfying the usual conditions.

If we have a map $f:P \sar Q$ of operads we immediately get a map
\begin{equation}
f^* : Q\text{-alg} \lar P\text{-alg}.
\end{equation}
If we introduce a simplicial direction, then we have a map $f_*$ in the other direction given by a simplicial bar construction. By that we mean that if $A$ is a $P$-algebra, then $B_\bullet(Q,P,A)$ is a $Q$-algebra, where $B_n(Q,P,A)=Q \circ P \circ \ldots \circ P \circ A$ with $P$ repeated $n$ times. If $\D$ has a notion of geometric realization we can get an honest $Q$-algebra $|B_\bullet(Q,P,A)|$ in $\D$. With some additional hypothesis it is possible to prove that $f_*$ is a homotopy left adjoint to $f^*$ and that if $f:P \sar Q$ is a weak equivalence of operads then $|B_\bullet(Q,P,A)|$ is weakly equivalent to $A$.

We say that an object $M$ in $\D$ is an $A$-bimodule if the sequence with $A$ in degree $0$, $M$ in degree $1$ and the initial object in all other degrees is a left $P$-module. Giving an $A$-bimodule structure on $M$ is equivalent to giving maps
\begin{equation} \label{PAmodule}
P(n) \otimes A^{\otimes i-1} \otimes M \otimes A^{\otimes n-i} \lar M
\end{equation}
for all $n \geq 1$ and $1 \leq i \leq n$ which satisfy the usual conditions.

The notions of a left or right $A$-module do not fit comfortably into this framework, but we will fix that in the next section.

\section{Subcategories of noncommutative sets} \label{subnoco}
Let $\noco$ be the category of noncommutative sets, as in \cite[section 6.1]{Loday} and \cite{PiRi}. The objects in $\noco$ are finite sets (the empty set is allowed) and a morphism is a map $f:S \sar T$ of finite sets together with a linear ordering of each $f^{-1}(t)$, $t \in T$. This category is noncommutative in the sense that there are $2$ different maps from a set with $2$ elements to a set with $1$ element, and as the notation suggests a map in $\noco$ can be factored uniquely as a permutation followed by a map in $\bD$ (if we pick a linear ordering of $S$ and $T$). This is exactly the data we need so that given an associative algebra $A$ and a morphism $f:S \sar T$ in $\noco$, we get a map $f_* : A^S \sar A^T$ in a natural way.

\subsection*{Enriching subcategories of noncommutative sets}
If we have a subcategory $\A$ of $\noco$, then each map $f:S \sar T$ in $\A$ comes with a linear ordering of each $f^{-1}(t)$. Thus, given an operad $P$ in $\C$, we can make sense of $P[f]=\bigotimes_{t \in T} P(f^{-1}(t))$ and we can define a new category $\A_P$ enriched over $\C$ as follows:

\begin{definition}
Let $\A$ be a subcategory of $\noco$ and let $P$ be an operad in $\C$. We define a category $\A_P$ enriched over $\C$ as follows. The objects are the same as in $\A$, and the $Hom$ objects are given by
\begin{equation}
Hom_{\A_P}(S,T)=\coprod_{f \in Hom_\A(S,T)} P[f].
\end{equation}
Composition is defined in terms of the operad structure in the evident way.
\end{definition}

For this definition to make sense, it is important that the monoidal product in $\C$ distributes over coproducts (see remark \ref{quasi_monoidal}). The composition is the map
\begin{equation}
\coprod_{g : T \sar U} P[g] \otimes \coprod_{f : S \sar T} P[f] \arrow{\cong} \coprod _{S \overset{f}{\sar} T \overset{g}{\sar} U} P[g] \otimes P[f] \lar \coprod_{h : S \sar U} P[h],
\end{equation}
where the first map uses that $\otimes$ distributes over $\coprod$ and the second map is given by the operad structure on $P$.

\begin{example}
There are several natural examples of subcategories of $\bD \Sigma$ one can consider; we list some of them here.
\begin{enumerate}
\item The simplicial category $\bD$, which we have defined to be the category of all nonempty finite totally ordered sets and order-preserving maps, or its augmentation $\bD_+$ which is the category $\bD$ together with the empty set.
\item The category $^0 \! \bD$ of finite totally ordered sets with a minimal element $0$ and order-preserving maps preserving the minimal element, or the category $^1 \! \bD$ of finite totally ordered sets with a maximal element $1$.
\item The doubly based category $^{01} \! \bD$ of finite totally ordered sets with both a minimal and maximal element.
\item Connes' cyclic category $\bD C$, which consists of finite cyclically ordered sets and order-preserving maps together with a linear ordering of each inverse image of an element, or its augmentation $\bD C_+$.
\item The category $^0 \! \bD C$ consisting of cyclically ordered sets with a basepoint $0$ and basepoint-preserving maps.
\end{enumerate}
\end{example}

\subsection*{Two generalizations of simplicial sets}
We start with the following easy, but important, lemma:
\begin{lemma}
The categories $^{01} \! \bD$ and $^0 \! \bD C$ are isomorphic to $\bD^{op}$.
\end{lemma}

\begin{proof}
(See e.g.~\cite[p.~ 621]{Dr04}.)
We construct a functor $^{01} \! \bD \sar \bD^{op}$ by sending a set with $n+2$ elements, say, $\{0,x_1,\ldots,x_n,1\}$ to $\mathbf{n}=\{0,1,\ldots,n\}$. The map on $Hom$ sets is given as in the following picture: 
\begin{equation}
\xymatrix{
0 \ar@{-}[rr]^-{0} \ar[rd] & & x_1 \ar@{-}[rr]^-{1} \ar[rd] & & x_2 \ar@{-}[rr]^-{2} \ar[ld] & & 1 \ar[ld] \\
& 0 \ar@{-}[rr]_-{0} & \ar[ul] & y_1 \ar@{-}[rr]_-{1} & \ar[ur] & 1 & 
}
\end{equation}

For the second case, given $f: S \sar T$ in $^0 \! \bD C$, the linear ordering of $f^{-1}(0)$ allows us to extend $f$ to a map $S \coprod \{1\} \sar T \coprod \{1\}$ by sending the elements in $f^{-1}(0)$ which are greater than $0$ to $0$ and the elements that are less than $0$ to $1$. The rest of the proof is the same as for the first case.
\end{proof}

Thus we get two generalizations of a simplicial object in $\D$, either as a functor $^{01} \! \bD_P \sar \D$ or a functor $^0 \! \bD C_P \sar \D$. As we will see, the first generalization allows us to define the $2$-sided bar construction while the second allows us to define the cyclic bar construction.

\begin{remark} \label{alloverC}
By a functor $F: \D \sar \D'$ between two categories that are enriched over $\C$, we will always mean a $\C$-functor, i.e., $F$ comes with maps $Hom_\D(X,Y) \sar Hom_{\D'}(F(X),F(Y))$ of objects in $\C$. For example, if $\C=\Top$ this means that we require all maps of $Hom$ sets to be continuous.
\end{remark}

If $X : {}^{01} \! \bD_P \sar \D$ or $^0 \! \bD C_P \sar \D$ we will abuse notation and denote the image of an element which corresponds to $\mathbf{n}=\{0,1,\ldots,n\}$ in $\bD$ by $X_n$. If $P(0)=P(1)=I$, the injective maps in $^{01} \! \bD$ or $^0 \! \bD C$, which correspond to surjective maps in $\bD$, do not change when we pass from $^{01} \! \bD$ to $^{01} \! \bD_P$ or from $^0 \! \bD C$ to $^0 \! \bD C_P$, and we will denote the map corresponding to $s_j$ in $\bD^{op}$ by $s_j$.

This gives a new way to look at a $P$-algebra $A$, as well as right, left and bimodules over $A$. A $P$-algebra $A$ gives a functor $(\bD_+)_P \sar \D$ given by $S \mapsto A^{\otimes S}$. A right $A$-module $M$ gives a functor $^0 \! \bD_P \sar \D$ given by $S \mapsto M \otimes A^{\otimes S-\{0\}}$, a left $A$-module $N$ gives a functor $^1 \! \bD_P \sar \D$ by $S \mapsto A^{\otimes S-\{1\}} \otimes N$, while a pair $(M,N)$ of a right and a left $A$-module gives a functor $^{01} \! \bD_P \sar \D$ by $S \mapsto M \otimes A^{\otimes S-\{0,1\}} \otimes N$. An $A$-bimodule $M$ gives a functor $^0 \! \bD C_P \sar \D$ by $S \mapsto M \otimes A^{\otimes S-\{0\}}$.

We will call the functor $S \mapsto M \otimes A^{\otimes S-\{0,1\}} \otimes N$ from $^{01} \! \bD_P$ to $\D$ the $P$-simplicial $2$-sided bar construction and denote it by $B^P_\bullet(M,A,N)$ or simply $B_\bullet(M,A,N)$ if $P$ is clear from the context. Similarly, we will call the functor $S \mapsto M \otimes A^{\otimes S-\{0\}}$ from $^0 \! \bD C_P$ to $\D$ the $P$-cyclic bar construction and denote it by $B^{cy}_\bullet(A;M)$. We will call the functor $S \mapsto Hom(A^{\otimes S-\{0\}},M)$ the cyclic cobar construction and denote it by $C_{cy}^\bullet(A;M)$. (Note that this terminology is not standard, as the usual cobar construction takes a coalgebra as input. We can think of the coalgebra in question as being the dual of $A$.)  This is a functor from $^0 \! \bD C_P^{op}$ to $\D$.

\begin{remark}
If $P$ is the associative operad in abelian groups, so a $P$-algebra is simply an algebra, then the homology of $B^{cy}_\bullet(A;M)$ is the Hochschild homology $HH_*(A;M)$ and the cohomology of $C_{cy}^\bullet(A;M)$ is the Hochschild cohomology $HH^*(A;M)$.
\end{remark}

\begin{notation}
We will call a functor $F$ from $\A_P$ to $\C$ or $\D$ a left $\A_P$-module and a functor $R$ from $\A_P^{op}$ to $\C$ or $\D$ a right $\A_P$-module.
\end{notation}

Thus a $P$-algebra $A$ gives rise to a left $\A_P$-module $S \mapsto A^{\otimes S}$ for any subcategory $\A$ of $\bD \Sigma$, and since a $P$-algebra is really a left $P$-module structure on a certain sequence under the composition product (Definition \ref{comp_prod}), a left $\A_P$-module is a generalization of this which depends on $\A$. (Though a general left $P$-module does not give a left $\A_P$-module in any natural way.)

\begin{remark}
The category $^0 \! \bD C_P$ should be compared to the category $\bar{\bD}^{op}$ and its enrichment $\hat{P}$ from \cite{Th79} (or $\widehat{\C}$, since $\C$ is an operad in Thomason's paper). Indeed, Thomason's category $\bar{\bD}^{op}$ is very closely related to the category $^0 \! \bD C$ of based cyclic sets. His condition in Definition 1.1 that for a map $f \in \bar{\bD}^{op}$, if $f(i_0)=0$ then either $f(i_1)>0$ only if $i_1>i_0$ or only if $i_1<i_0$ guarantees that $f$ can be lifted to a map of cyclically ordered sets. But he only has one map from $\mathbf{n}=\{0,1,\ldots,n\}$ to $\mathbf{0}=\{0\}$ instead of $n+1$ maps, so the lift is not always unique. Moreover, his enriched category $\hat{P}$ is the same as our $\bD C_P$, except that he only uses the spaces $P(f^{-1}(t))$ for $t \not = 0$. These differences do not matter as long as one only studies the case when $M=I$, which is in effect what he did.
\end{remark}

We also note that a functor $R : (\bD_+)_P^{op} \sar \C$ is precisely a right $P$-module. Thus we also obtain various generalizations of right $P$-modules. For example, we will think of a right $(\bD C_+)_P$-module as a right $P$-module with some extra structure.

\subsection*{Geometric realization}
We start with some notation.

\begin{definition} \label{geom_real_def}
If we have a left $\A_P$-module $F$ and a right $\A_P$-module $R$, we define $R \otimes_{\A_P} F$ as the coequalizer
\begin{equation}
\coprod_{[f:S \rightarrow T]} R(T) \otimes P[f] \otimes F(S) \rightrightarrows \coprod_{[S]} R(S) \otimes F(S) \lar R \otimes_{\A_P} F.
\end{equation}

Similarly, if $F$ and $R$ are both right $\A_P$-modules we define $Hom_{\A_P}(R,F)$ as the equalizer
\begin{equation}
Hom_{\A_P}(R,F) \lar \coprod_{[S]} Hom(R(S),F(S)) \rightrightarrows \coprod_{[f:S \rightarrow T]} Hom(R(T) \otimes P[f],F(S)).
\end{equation}
\end{definition}

If $P=Ass$ and $\A$ is either $^{01} \! \bD$ or $^0 \! \bD C$, so that a left $\A_P$-module is a simplicial object, then $\Delta^\bullet \otimes_{\A_P} X_\bullet = |X_\bullet|$, where $\Delta^\bullet(S)$ is either $\Delta^{|S-\{0,1\}|}$ or $\Delta^{|S-\{0\}|}$.

If we have a map $f:P \sar Q$ of operads, then we get a functor $f: \A_P \sar \A_Q$ between enriched categories. Thus we can pull a functor $\A_Q \sar \D$ or $\A_Q^{op} \sar \D$ back to a functor $\A_P \sar \D$ or $\A_P^{op} \sar \D$. This gives a functor $f^* : \D^{\A_Q} \sar \D^{\A_P}$. In particular, since $\A_{Ass} \cong \A$, given any $A_\infty$ operad $P$ and a functor $\A \sar \D$ or $\A^{op} \sar \D$ we can pull it back to a functor $\A_P \sar \D$ or $\A_P^{op} \sar \D$.

With $\A={}^{01} \! \bD_P$ or $^0 \! \bD C_P$ we see that we can regard any simplicial object in $\D$ as a left $^{01} \! \bD_P$ or $^0 \! \bD C_P$-module. Similarly, any cosimplicial object can be regarded as a right $^{01} \! \bD_P$ or $^0 \! \bD C_P$-module.

\section{The associahedra and cyclohedra} \label{asscyclo}
The original definition of the associahedra can be found in \cite{St63}. The cyclohedra got their name from Stasheff \cite{St95}, but had been considered earlier, first by Bott and Taubes in \cite{BoTo94}. They sometimes go under the name Stasheff associahedra of type $B$. For an introduction to the associahedra and cyclohedra, see \cite{Markl}.

\subsection*{The associahedra}
\begin{figure}[h]
\begin{center}
\epsfig{file=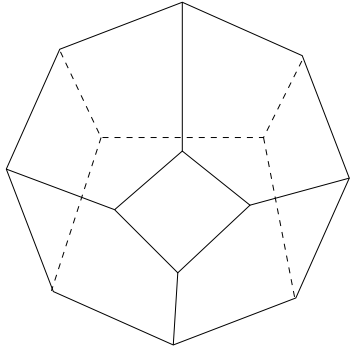, scale=0.6}
\caption{The associahedron $K_5$.}\label{K5}
\end{center}
\end{figure}

Consider all ways to parenthesize (in a meaningful way) $n$ linearly ordered variables. The maximal number of pairs of parentheses is $n-2$, and the Stasheff associahedron $K_n$ has an $(n-2-i)$-cell for each way to parenthesize using $i$ pairs of parentheses, with a face map for each way to insert an additional pair. Perhaps a more precise definition of $K_n$ is as the cone on $L_n$, where $L_n$ is the union of various copies of $K_r \times K_{n-r+1}$, as in \cite[Definition 1.7]{BoVo73}. This definition also makes sense in the category of simplicial sets. For example, $K_0=K_1=K_2=*$, $K_3=I$ is an interval and $K_4$ is a pentagon. Figure \ref{K5} shows $K_5$.

\begin{theorem} (Stasheff \cite{St63})
For $n \geq 2$ the associahedron $K_n$ is homeomorphic to $D^{n-2}$, and the sequence $\K=\{K_n\}_{n \geq 0}$ forms an $A_\infty$ operad in $\Top$.
\end{theorem}

The operad structure can be thought of as substitution of parenthesized expressions.

The Boardman-Vogt $W$-construction, which is usually given in terms of certain metric trees, provides a cofibrant replacement of operads in a certain model category. The associahedra operad $\K$ (or a cubical decomposition of it) is the $W$-construction on the associative operad $Ass$ \cite[Example 2.22]{MSS}, though we have to be careful about exactly which category of operads we work in.

The model category structure on operads (in $\Top$) is given by levelwise weak equivalences and fibrations, while the cofibrations are what they have to be \cite[Theorem 3.2]{BeMo}. This relies on special properties of the category $\Top$, though it is possible to weaken the conditions necessary for getting a model category somewhat by considering reduced operads \cite[Theorem 3.1]{BeMo}. An operad $P$ is reduced if $P(0)=I$. It is not clear if there is a model category structure on operads in a general symmetric monoidal model category.

But $\K$ is not cofibrant in the category of operads in $\Top$. For example, it is easy to see that there can be no map from $\K$ to the little intervals operad $\C_1$. The problem is that $* \in \C_1(0)$ does not act as a unit. If we perform the $W$-construction on $Ass$ in this category we also get an operad with a very big first space. The solution in \cite{MSS} is to consider the category of operads with $P(0)=\varnothing$. In this category one can show that there is a map from the associahedra to the little intervals operad, and the associahedra operad is indeed cofibrant in this category.

Because we want to generalize the $W$-construction in such a way that it produces the cyclohedra, we will sketch the details of the $W$-construction on $Ass$ in the topological setting. To do this we need to discuss trees. We will use \cite[Definition 3.1]{Ching05} as our definition of a tree. In particular, the root has exactly one incoming edge, and by a vertex we mean an internal vertex, i.e., we do not include the root or any of the leaves in the list of vertices. Because we are not considering $\Sigma$-operads, we also require that the leaves come with a linear ordering.

\begin{remark}
If we want to work in the category of operads where $P(0)$ is allowed to be nonempty, we need to allow trees where a vertex can have zero incoming edges. With the restriction that $P(0)=\varnothing$, we can restrict our attention to trees where each vertex has at least one incoming edge. This cuts the number of trees with a fixed number of leaves down dramatically.
\end{remark}

A metric on a tree $T$ is an assignment of a length $0 \leq w(e) \leq 1$ of each internal edge $e$ in $T$. Let $\mathcal{T}_n$ be the space of metric trees with $n$ leaves. Given a tree $T \in \mathcal{T}_n$ let $vert(T)$ be the set of internal vertices, and for $v \in vert(T)$ let $In(v)$ be the set of incoming edges to $v$.

\begin{figure}[h]
\begin{center}
\epsfig{file=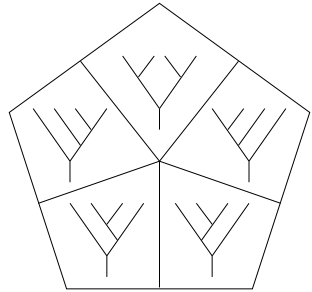, scale=0.7}
\caption{The associahedron $K_4$ from the $W$-construction.} \label{subK4}
\end{center}
\end{figure}

\begin{definition}
Given an operad $P$ in $\Top$ with $P(0)=\varnothing$, $WP$ is the operad defined as follows. $WP(n)$ is a quotient of
\begin{equation}
\coprod_{T \in \mathcal{T}_n} \prod_{v \in vert(T)} P(In(v))
\end{equation}
under the following identifications:
\begin{enumerate}
\item If an internal edge has length $0$, the tree is identified with the one obtained by contracting this edge and applying the corresponding operation in $P$.
\item If a vertex $v$ with one incoming edge is labeled by the unit $* \in \Top$, the tree is identified with the one obtained by deleting $v$ and giving the new edge the maximum length of the incoming and outgoing edge of $v$.
\end{enumerate}
The structure maps in $WP$ are given by grafting trees, and assigning the length $1$ to any new internal edges.
\end{definition}
For a much more general approach to the $W$-construction, see \cite{BeMo2}. The claim in \cite{MSS} is that the $W$-construction on $Ass$ in this setting gives $\K$. This is easy to verify, after noting that because $Ass(1)=*$ we can disregard trees where a vertex has only one incoming edge. For example, $K_4$ is given by Figure \ref{subK4}.

It is easy to see that $\K$ is cofibrant in the category of operads without zeroth space. This amounts to showing that for any trivial fibration $P \arrow{\simeq} Ass$ of operads, which just amounts to requiring that each $P(n)$ is fibrant and contractible, the dotted arrow in the diagram
\begin{equation} \xymatrix{
 & P \ar[d]^\simeq \\
{\K} \ar@{.>}[ur] \ar[r] & Ass}
\end{equation}
exists. But we can construct a map like this by induction. If we are given maps $K_i \sar P(i)$ for $i <n$, the map $K_n \sar P(n)$ is determined on $\partial K_n$, and now we just have to solve the extension problem
\begin{equation}
\xymatrix{
S^{n-3} \ar[r] \ar[d] & P(n) \ar[d] \\
D^{n-2} \ar@{.>}[ur] \ar[r] & {*}}
\end{equation}
which we know we can solve because each $P(n)$ is fibrant and contractible.

Thus for any $A_\infty$ operad $P$, we get a map $\K \sar P$ of operads without a zeroth space. If $P(0)=*$ and $P(0)$ acts as a unit in the sense that all diagrams of the form
\begin{equation}
\xymatrix{ K_n \times K_0 \ar[r] \ar[d]_{\circ_i} & P(n) \times P(0) \ar[d]^{\circ_i} \\
K_{n-1} \ar[r] & P(n-1)}
\end{equation}
commute, we can promote $\K \sar P$ to a map of operads with zeroth spaces, and then pull a functor $F: \A_P \sar \Top$ back to a functor $\A_\K \sar \Top$. Thus it makes sense to concentrate on the $A_\infty$ operad $\K$ as long as we are willing to restrict the class of operads we consider. This restriction excludes operads like the little intervals operad, so for some purposes this restriction is bad, for example if we want to consider tensor products of operads.

The operad $\K$ has an obvious filtration, where we let $\K_n$ be the operad generated by $K_i$ for $i \leq n$. An algebra over $\K_n$ is precisely an $A_n$ algebra as in \cite{St63}, and giving an $A_n$ algebra structure on $A$ is equivalent to giving maps $K_i \times A^i \sar A$ for $0 \leq i \leq n$ satisfying the usual conditions.

\subsection*{Associahedra as configuration spaces}
We can interpret the associahedra in terms of configuration spaces. Observe that the standard $n$-simplex $\Delta^n$ is the configuration space of $n+2$ points on the unit interval, where the first point is at $0$ and the last at $1$. Let the points be labeled by elements of a set $S \in {^{01} \! \bD}$ with $|S|=n+2$, say, $S=\{0,x_1,\ldots,x_n,1\}$. By abuse of notation, we will let $x_i$ denote both the position of the $(i+1)$st point and an element in $S$. We get the standard description of $\Delta^n$ by setting $t_i=x_{i+1}-x_{i}$ ($t_0=x_1$, $t_n=1-x_n$). The associahedron $K_{n+2}=\K(S)$ is also such a configuration space; it is the compactification of the configuration space of $n+2$ distinct points as above, where for each time we have a point repeated $i$ times we use a copy of $K_i$ instead of just a point. This gives an inductive definition of the associahedra. This is the Axelrod-Singer compactification, see \cite{AxSi} and also \cite{FuMa}. Of course, we could also describe the configuration space of distinct points on $I$ in terms of distinct points on $\R$ modulo translation and dilatation. This is the point of view found in the references, and the point of view we have to take if we want to generalize to configurations on higher-dimensional manifolds, see Remark \ref{higherconf}.

If $f:S \sar T$ is a map in $^{01} \! \bD$, we can interpret $\K(T) \otimes \K[f] \sar \K(S)$ in terms of the above configuration space as follows. For each $t \in T$, the map replaces the point labeled by $t$ by points labeled by the set $f^{-1}(t)$, and the factor $\K(f^{-1}(t))$ tells us how. This works because $0 \in f^{-1}(0)$ and $1 \in f^{-1}(1)$, so we never remove the points at the beginning and end of the interval. 

Next we compare the associahedra to the standard $n$-simplexes. Denote by $s^i : K_{n+2} \sar K_{n+1}$, $0 \leq i \leq n-1$, the map $K_{n+2} \cong K_{n+2} \times K_0 \sar K_{n+1}$ obtained from $\circ_{x_{i+1}}$ (Definition \ref{circle_notation}). This is the same as removing the $(i+2)$nd variable in the parenthesized expression of $n+2$ variables defining $K_{n+2}$, and as we just saw it corresponds to removing the point marked $x_{i+1}$ in the configuration space. We also get maps $K_{n+2} \sar K_{n+1}$ by removing $0$ or $1$, but these do not correspond to codegeneracy maps on $\Delta^n$. Similarly we get maps $d^j : K_{n+2} \cong K_{n+2} \times K_2 \sar K_{n+3}$ for $0 \leq j \leq n+1$ from $\circ_{x_j}$ ($x_0=0$, $x_{n+1}=1$), which correspond to replacing $x_j$ with a double point. 

As is obvious from the configuration space interpretation of $K_n$, there is a surjective map $K_{n+2} \sar \Delta^n$ which is a homeomorphism on the interior. The association $n \mapsto K_{n+2}$ is not quite a cosimplicial space, because some of the simplicial identities commute only up to homotopy, but the following diagrams commute:
\begin{equation} \label{mapsKnDeltan} 
\xymatrix{
K_{n+2} \ar[r] \ar[d]^{s^i} & \Delta^n \ar[d]^{s^i} & & K_{n+2} \ar[r] \ar[d]^{d^j} & \Delta^n \ar[d]^{d^j} \\
K_{n+1} \ar[r] & \Delta^{n-1} & & K_{n+3} \ar[r] & \Delta^{n+1} } 
\end{equation}

In particular these two diagrams show that $S \mapsto \Delta^{|S|-1}$ gives a functor $^{01}\!\bD_\K^{op} \sar \Top$. We also see that the two faces of $K_{n+2}$ coming from the inclusions $K_2 \times K_{n+1} \sar K_{n+2}$ are crushed to a point in $\Delta^n$.

Let $A$ be a $\K$-algebra, $M$ a right $A$-module and $N$ a left $A$-module. By regarding $\K$ as a functor $^{01}\! \bD_\K^{op} \sar \Top$ we can now define the $2$-sided bar construction as
\begin{equation}
B(M,A,N)= \K \otimes_{{^{01}}\!\bD_\K} B^\K_\bullet(M,A,N).
\end{equation}
We elaborate on what this means. The tensor product is defined as a coequalizer, which just means a quotient in $\Top$, so $B(M,A,N)$ is given by
\begin{equation}
\coprod_n K_{n+2} \times M \times A^n \times N/\sim,
\end{equation}
where we identify $(f^* x,y)$ with $(x,f_* y)$. (The roles of $f_*$ and $f^*$ are reversed here because $^{01} \! \bD_\K$ has replaced $\bD^{op}$, not $\bD$.) 

\begin{proposition}
When $M=N=*$, the bar construction as defined above agrees with Stasheff's definition in \cite{St63}. In particular, if $A$ is grouplike ($\pi_0 A$ is a group) then $BA=B(*,A,*)$ is a delooping of $A$.
\end{proposition}

Similarly, given any functor $X : {^{01}} \! \bD_\K \sar \Top$ we define $|X|$ as $\K \otimes_{{^{01}} \! \bD_\K} X$, and given a functor $Y : {^{01}} \! \bD_\K^{op} \sar \Top$ we define $Tot(Y)$ as $Hom_{{^{01}} \! \bD_\K}(\K,Y)$.

\begin{proposition} \label{KDeltanatural}
The maps $K_{n+2} \sar \bD^n$ assemble to a natural transformation of functors from $^{01} \! \bD_\K$ to $\Top$.
\end{proposition}

\begin{proof}
This follows immediately from the two commutative diagrams in (\ref{mapsKnDeltan}).
\end{proof}

\begin{remark} \label{rem:representability}
Consider the functor
\begin{equation}
Hom_{{}^{01} \! \bD_\K}(\{0,x_1,\ldots,x_n,1\},-) : {}^{01} \! \bD_\K \lar \Top.
\end{equation}
It has ``$0$-simplices'' $\coprod_{1 \leq i \leq n+1} K_i \times K_{n+2-i}$ etc, up to exactly one nondegenerate ``$n$-simplex''. But this is a description of $K_{n+2}$, and explains why we need to use $K_{n+2}$ instead of $\Delta^n$ to get a good notion of geometric realization.
\end{remark}

\begin{proposition} \label{comparinggeomreal}
Let $X$ be a simplicial space and regard $X$ as a functor $^{01} \! \bD_\K \sar \Top$ via
\begin{equation}
^{01} \! \bD_\K \lar {^{01}} \! \bD_{Ass} \arrow{\cong} \bD^{op} \arrow{X} \Top.
\end{equation}
Then the natural map
\begin{equation}
\K \otimes_{{^{01}} \! \bD_\K} X \lar \Delta^\bullet \otimes_{\bD^{op}} X
\end{equation}
is a homeomorphism.

Similarly, if we regard a cosimplicial space $Y$ as a functor $^{01} \! \bD_\K^{op} \sar \Top$ the natural map
\begin{equation}
Hom_{\bD}(\Delta^\bullet,Y) \lar Hom_{{}^{01} \! \bD_\K}(\K,Y)
\end{equation}
is a homeomorphism.
\end{proposition}

\begin{proof}
Let $S=\{0,x_1,\ldots,x_n,1\}$, $T=\{0,y_1,\ldots,y_m,1\}$ and let $f:S \sar T$ be a map in $^{01}\! \bD$ which is dual to $g:\mathbf{m} \sar \mathbf{n}$ in $\bD$. Then the component of $f$ in $^{01} \! \bD_\K$ is $\K[f]$, and $\bar{f}_*=g^* : X_n \sar X_m$ for any $\bar{f} \in Hom_{{}^{01} \! \bD_\K}(S,T)$ in the component of $f$. Thus the induced map
\begin{equation}
\K[f] \times X_n \lar X_m
\end{equation}
is constant in the first variable, and when we consider the corresponding map $\K(T) \times \K[f] \sar \K(S)$ which appears in the coequalizer defining $\K \otimes_{{^{01}}\! \bD_\K} X$, the image of $\{k\} \times K[f]$ in $K_{n+2} \times X_n$ is crushed to a point. Doing this for all $f$ gives us exactly the projections $K_{n+2} \sar \Delta^n$.

The second part is similar.
\end{proof}

\begin{remark} \label{An_R_remark}
We would like to know that under reasonable circumstances a map $X \sar Y$ of functors from $^{01} \! \bD_\K$ to $\Top$ which is a levelwise weak equivalence gives an equivalence after geometric realization, or that the skeletal filtration of $X$ gives rise to a spectral sequence. We will return to this point in \cite{An_R}, where we define an enriched Reedy category and show that $^{01} \! \bD_\K$ is an example. In this example a levelwise weak equivalence between Reedy cofibrant functors gives a weak equivalence after geometric realization, and the skeletal filtration of a Reedy cofibrant functor gives the expected spectral sequence.
\end{remark}

\subsection*{The cyclohedra}
\begin{figure}[h]
\begin{center}
\epsfig{file=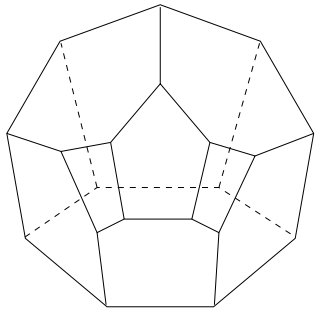, scale=0.6}
\caption{The cyclohedron $W_4$.} \label{W4}
\end{center}
\end{figure}

Next we consider the right module over $\K$ given by the cyclohedra. To define the cyclohedron $W_n$, we again consider all ways to parenthesize $n$ variables, but now we let them be cyclically ordered. In this case the maximal number of pairs of parentheses is $n-1$. For example, $12$ can be parenthesized as either $(12)$ or $(21)$. The same construction as above, now with an $n-i-1$ cell for each way to parenthesize $n$ variables using $i$ pairs of parentheses gives the space $W_n$. Alternatively we can define $W_n$ as the cone on a union of various copies of $W_s \times K_{n-s+1}$. For example, $W_0=W_1=*$, $W_2=I$ and $W_3$ is a hexagon. Figure \ref{W4} shows $W_4$. 

\begin{theorem} \label{Wiscyclic}
For $n \geq 1$ the cyclohedron $W_n$ is homeomorphic to $D^{n-1}$, and the cyclohedra assemble to a functor $\mathcal{W} : (\bD C_+)_\K^{op} \sar \Top$.
\end{theorem}

\begin{proof}
It is well known (\cite[\S 4]{St95}) that $W_n \cong D^{n-1}$ for $n \geq 1$ and that $\mathcal{W}$ is a right $\K$-module, i.e., a functor $(\bD_+)_\K^{op} \sar \Top$.

To show that $\mathcal{W}$ is a functor from $(\bD C_+)_\K^{op}$ to $\Top$, we need to give a map
\begin{equation}
\mathcal{W}(T) \otimes \K[f] \lar \mathcal{W}(S)
\end{equation}
for each map $f : S \sar T$ in $\bD C_+$. But a map $f : S \sar T$ corresponds to a partition of $S$ into sets $f^{-1}(t)$ for $t \in T$ (some of these may be empty). Such a partition corresponds to a way to parenthesize the elements in $S$, by putting a parenthesis around $f^{-1}(t)$ if $|f^{-1}(t)| \geq 2$. This in turn corresponds to a cell in $\mathcal{W}(S)$, and the map $\mathcal{W}(T) \otimes \K[f] \sar \mathcal{W}(S)$ is the inclusion of this cell.
\end{proof}

\begin{remark}
Theorem \ref{Wiscyclic} is related to Markl's result \cite[Theorem 2.12]{Markl} that if we consider the $\Sigma$-operad with $n \mapsto \Sigma_n \times K_n$ then the symmetric sequence $n \mapsto C_n \backslash \Sigma_n \times W_n$ is a right module over this operad.
\end{remark}

Next we describe the version of the Boardman-Vogt $W$-construction which gives us $\mathcal{W}$. The idea is to change the definition of a tree slightly, in a way that is similar to \cite[Definition 7.3]{Ching05}. We require the leaves to come with a cyclic ordering rather than a linear ordering, and we allow the root to have more than one incoming edge. A tree $T$ comes with a partition of the leaves according to which connected component of $T-\{r\}$ they belong to, and part of the data of a tree is a compatible linear ordering of the leaves belonging to each of the root edges.

We also regard the root edges as internal edges, and assign a length to each of them. Let $\mathcal{T} C_n$ be the space of such trees with leaves labeled by the cyclically ordered set $\{1,2,\ldots,n\}$. Note that the open subset of $\mathcal{T} C_n$ where the root has exactly one incoming edge is dense in $\mathcal{T} C_n$, and this subspace is the disjoint union of $n$ copies of $\mathcal{T}_n \times I$. Also let $\overline{\bD C}_P$ be the enrichment of the category of cyclically ordered sets and order-preserving \emph{surjections}.

\begin{definition}
Given an operad $P$ in $\Top$ with $P(0)=\varnothing$ and a functor $R : \overline{\bD C}^{op}_P \sar \Top$, $W(P;R)$ is the functor from $\overline{\bD C}^{op}_{WP}$ to $\Top$ defined as follows. The space $W(P;R)(n)$ is a quotient of
\begin{equation}
\coprod_{T \in \mathcal{T} C_n} R(In(r)) \times \Big( \prod_{v \in vert(T)} P(In(v)) \Big)
\end{equation}
under the following identifications:
\begin{enumerate}
\item If an internal edge has length $0$, the tree is identified with the one obtained by contracting this edge and applying the corresponding operation from the operad structure on $P$ or the right $\overline{\bD C}_P$-module structure on $R$.
\item If a vertex $v \neq r$ with one incoming edge is labeled by the unit $* \in \Top$, the tree is identified with the one obtained by deleting $v$ and giving the new edge the maximum length of the incoming and outgoing edge of $v$.
\end{enumerate}
The structure maps are given by grafting trees, and assigning the length $1$ to any new internal edges.
\end{definition}

\begin{proposition}
The cyclohedra can be obtained as $\mathcal{W}=W(Ass; R)$, where $R : \overline{\bD C}^{op}_{Ass} \sar \Top$ sends any set to $*$.
\end{proposition}

\begin{figure}
\begin{center}
\epsfig{file=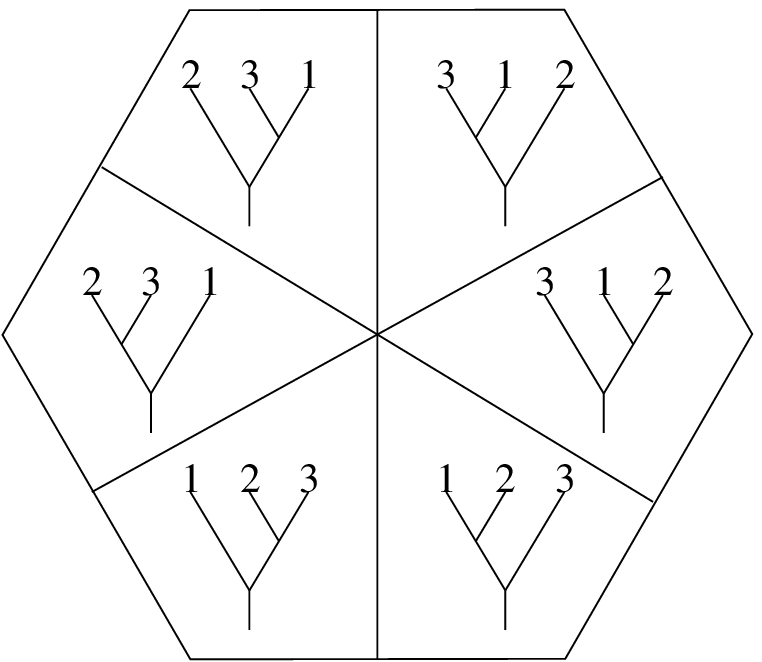, scale=0.69}
\caption{The cyclohedron $W_3$ from a relative $W$-construction.} \label{subcyclo3}
\end{center}
\end{figure}

\begin{proof}
The relative $W$-construction gives an $(n-1)$-cube for each binary tree with only one incoming root edge. We can think of each vertex of $W_n$ as a binary tree where all the internal edges have length $1$. By subdividing $W_n$ as in Figure \ref{subcyclo3} we see that $W_n$ can be decomposed as a union of $(n-1)$-cubes, one for each binary tree.
\end{proof}

\begin{figure}[h]
\begin{center}
\epsfig{file=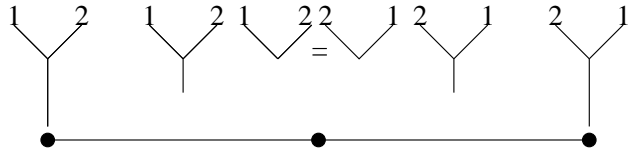, scale=0.9}
\caption{The cyclohedron $W_2$ from a relative $W$-construction.}
\end{center}
\end{figure}

The cyclohedra, regarded as a functor $\overline{\bD C}^{op}_\K \sar \Top$, also has a universal lifting property. For any $R : \overline{\bD C}^{op}_\K \sar \Top$ with each $R(S) \simeq *$, we can construct a natural transformation $\mathcal{W} \sar R$ by lifting one cell at a time.

\begin{remark}
We could consider this construction for other subcategories of $\bD \Sigma$. For example, if we consider $R : \bar{\bD \Sigma}^{op}_{Ass} \sar \Top$ given by $R(S)=*$ for all $S$ we find that $W(Ass;R)(3)$ is two copies of $W_3$ joined at the center. More generally, $W(Ass;R)(n)$ is the union of $(n-1)!$ copies of $W_n$ along a codimension $2$ subspace.
\end{remark}

The functor $\mathcal{W}$ has a filtration which is compatible with the filtration of $\K$. We let $\mathcal{W}_n$ be the functor $\mathcal{W}_n : (\bD C_+)_{\K_n}^{op} \sar \Top$ generated $W_i$ for $i \leq n$. 

Again we can relate this to configuration spaces. We can also consider the $n$-simplex $\Delta^n$ as the configuration space of $n+1$ points on $S^1$ labeled by elements of some $S \in {^0}\! \bD C$ with $|S|=n+1$, say, $S=\{0,x_1,\ldots,x_n\}$, with $0$ at the basepoint. Then $W_{n+1}$ is the Axelrod-Singer compactification of the interior of this space where we use a copy of $K_i$ instead of a point every time we have a point repeated $i$ times. 

Again we have maps
\begin{equation}
s^i : W_{n+1} \lar W_n
\end{equation}
for $0 \leq i \leq n-1$, which can be interpreted as removing the point $x_{i+1}$ in the above configuration space, and maps $d^j : W_{n+1} \sar W_{n+2}$ for $0 \leq j \leq n+1$ which can be interpreted as replacing the point $x_i$ with a double point. Here $d^0$ and $d^{n+1}$ both replace $0$ by a double point, either with $\{0,x_1\}$ or $\{x_{n+1},0\}$. Again we have a surjective map $W_{n+1} \sar \Delta^n$, which gives commutative diagrams
\begin{equation} 
\xymatrix{
W_{n+1} \ar[r] \ar[d]^{s^i} & \Delta^n \ar[d]^{s^i} & & W_{n+1} \ar[r] \ar[d]^{d^j} & \Delta^n \ar[d]^{d^j} \\
W_n \ar[r] & \Delta^{n-1} & & W_{n+2} \ar[r] & \Delta^{n+1} } 
\end{equation}

We also get a map $W_{n+1} \sar W_n$ from removing the point $0$, but this does not correspond to a codegeneracy map on $\Delta^n$.

If we have a functor $X : {^0}\! \bD C_\K \sar \Top$, we define the geometric realization $|X|$ as $\mathcal{W} \otimes_{{^0}\! \bD C_\K} X$, and if we have a functor $Y: {}^0\! \bD C_\K^{op} \sar \Top$ we define the total space $Tot(Y)$ as $Hom_{{}^0 \! \bD C_\K}(\mathcal{W},Y)$. In particular, given an $A_\infty$ $H$-space $A$ and an $A$-bimodule $M$ we can define the cyclic bar construction $B^{cy}(A;M)$ and the cyclic cobar construction $C_{cy}(A;M)$ this way.

It is clear that the analogs of Proposition \ref{KDeltanatural} and \ref{comparinggeomreal} hold:

\begin{proposition} \label{WDeltanatural}
The maps $W_{n+1} \sar \bD^n$ assemble to a natural transformation of functors from $^0 \! \bD C_\K$ to $\Top$.
\end{proposition}

\begin{proposition} \label{comparinggeomreal2}
Let $X$ be a simplicial space and regard $X$ as a functor $^0 \! \bD C_\K \sar \Top$ via
\begin{equation}
^0 \! \bD C_\K \lar {}^0 \! \bD C_{Ass} \arrow{\cong} \bD^{op} \arrow{X} \Top.
\end{equation}
Then the natural map
\begin{equation}
\mathcal{W} \otimes_{{}^0 \! \bD C_\K} X \lar \Delta^\bullet \otimes_{\bD^{op}} X
\end{equation}
is a homeomorphism.

Similarly, if we regard a cosimplicial space $Y$ as a functor $^0 \! \bD C_\K^{op} \sar \Top$ the natural map
\begin{equation}
Hom_{\bD}(\Delta^\bullet,Y) \lar Hom_{{}^0 \! \bD C_\K}(\mathcal{W},Y)
\end{equation}
is an homeomorphism.
\end{proposition}

\begin{remark}
As in Remark \ref{rem:representability}, we can consider the functor
\begin{equation}
Hom_{{}^0 \! \bD C_\K}(\{0,x_1,\ldots,x_n\},-) : {}^0 \! \bD C_\K \lar \Top.
\end{equation}
It has ``$0$-simplices'' $\coprod_{n+1} K_{n+1}$ etc, up to exactly one nondegenerate ``$n$-simplex''. This time we get $W_{n+1}$, and this explains why we need to use $W_{n+1}$ instead of $\Delta^n$ (or $K_{n+2}$) to get a good notion of geometric realization.
\end{remark}

We can also consider functors defined only on $^0 \! \bD C_{\K_n}$. If $X: {}^0 \! \bD C_{\K_n} \sar \Top$, the expression $\mathcal{W}_n \otimes_{{}^0 \! \bD C_{\K_n}} X$ makes sense, and by abuse of notation we will denote it by $sk_{n-1} |X|$, because if $X$ is the restriction of a functor from $^0 \! \bD C_\K$ then this does give the $(n-1)$-skeleton. Similarly we will denote $Hom_{{}^0 \! \bD C_{\K_n}}(\mathcal{W}_n,Y)$ by $Tot^{n-1}(Y)$ for a functor $Y : {}^0 \! \bD C_{\K_n}^{op} \sar \Top$. In particular, given a pair $(A,M)$ consisting of an $A_n$ algebra $A$ and an $A$-bimodule $M$ we can define $sk_{n-1} B^{cy}(A;M)$ and $Tot^{n-1} C_{cy}(A;M)$.

\begin{remark}
In \cite{An_R} we will also show that $^0 \! \bD C_\K$ is an enriched Reedy category and that the same things we mentioned in Remark \ref{An_R_remark} hold in this case.
\end{remark}

\begin{remark} \label{higherconf}
A natural generalization is to consider the configuration space of points in $\R^m$ modulo translation and dilatation, and the resulting $\Sigma$-operad $F_m$ we get from the Axelrod-Singer compactification of this space. If we have a parallelizable $m$-manifold $N$, the compactified configuration space of points on $N$ is naturally a right module over $F_m$. If we want to consider a right $F_m$-algebra $A$ and an $(F_m,A)$-module $M$, it is natural to consider a pointed manifold and based configurations. If $N$ is not parallelizable one has to consider a framed version of this, we refer to \cite{Markl_comp} for the details.
\end{remark}

\section{Traces and cotraces} \label{tracecotrace}
In this last section we identify what kind of structure the cyclic bar or cobar construction corepresents or represents. This perspective will be convenient in our upcoming paper \cite{AnTHH}.

\subsection*{Traces on a $P$-algebra $A$}
We start by recalling Markl's definition of a trace (\cite[Definition 2.6]{Markl}), adapted to right modules over operads (as opposed to $\Sigma$-operads). As usual, let $P$ be an operad and let $A$ be a $P$-algebra. Let $\E_A$ be the endomorphism operad for $A$, with $\E_A(S)=Hom(A^{\otimes S},A)$, and let $\E_{A,B}$ be the sequence given by $\E_{A,B}(S)=Hom(A^{\otimes S},B)$. Then $\E_{A,B}$ is a right $\E_A$-module, and by using the map $P \sar \E_A$ defining the $P$-algebra structure on $A$, a right $P$-module. Given another right $P$-module $R$, we can ask for a map $R \sar \E_{A,B}$ of right $P$-modules. Markl defines an $R$-trace on $A$ (with target $B$) as such a map. This is equivalent to giving maps $R(S) \otimes A^{\otimes S} \sar B$ for each finite totally ordered set $S$, such that the diagram
\begin{equation}
\xymatrix{
{R}(T) \otimes P[f] \otimes A^{\otimes S} \ar[r] \ar[d] & {R}(T) \otimes A^{\otimes T} \ar[d] \\
{R}(S) \otimes A^{\otimes S} \ar[r] & B }
\end{equation}
commutes for all maps $f:S \sar T$ in $\bD_+$. Thus an $R$-trace is corepresented by $R \otimes_{(\bD_+)_P} A^\bullet$, where $A^\bullet$ is the functor $S \mapsto A^{\otimes S}$. In other words, an $R$-trace on $A$ with target $B$ is the same thing as a map 
\begin{equation}
R \otimes_{(\bD_+)_P} A^\bullet \lar B.
\end{equation}

\subsection*{Traces on a pair $(A,M)$}
We will modify this construction so that it applies to the situation where we have a pair $(A,M)$ consisting of a $P$-algebra $A$ and an $A$-bimodule $M$. Instead of a right $P$-module we now need a functor  $R : {^0 \! \bD C_P^{op}} \sar \C$. Let $S \in {^0 \! \bD C}$, and let $\E_{A,M,B}$ be the functor $^0 \! \bD C_P^{op} \sar \C$ defined by $\E_{A,M,B}(S)=Hom(M \otimes A^{\otimes S-\{0\}},B)$. Then we can ask for a natural transformation $R \sar \E_{A,M,B}$ of functors from $^0 \! \bD C_P^{op}$ to $\C$.

\begin{definition}
Let $A$ be a $P$-algebra, $M$ an $A$-bimodule, and let $R$ be a functor ${^0 \! \bD C_P^{op}} \sar \C$. An $R$-trace on $(A,M)$ with target $B$ is a natural transformation $R \sar \E_{A,M,B}$ of functors.
\end{definition}

Giving an $R$-trace is equivalent to giving maps $R(S) \otimes M \otimes A^{\otimes S-\{0\}} \sar B$ for each $S \in {^0 \! \bD C}$ such that the diagram
\begin{equation} \label{trace_comm}
\xymatrix{
{R}(T) \otimes P[f] \otimes M \otimes A^{\otimes S-\{0\}} \ar[r] \ar[d] &  {R}(T) \otimes M \otimes  A^{\otimes T-\{0\}} \ar[d] \\
{R}(S) \otimes M \otimes  A^{\otimes S-\{0\}} \ar[r] & B }
\end{equation}
commutes for all $f:S \sar T$ in $^0 \! \bD C$.

As in the classical case, it is easy to find the object corepresenting $R$-traces. Indeed, it follows immediately from the definitions that giving an $R$-trace on $(A,M)$ with target $B$ is equivalent to giving a map
\begin{equation}
R \otimes_{^0 \! \bD C_P} B^{cy}_\bullet(A;M) \lar B.
\end{equation}

In particular, we have the following:

\begin{observation}
With $P=\K_n$, $1 \leq n \leq \infty$, a $\mathcal{W}_n$-trace is corepresented by the partial cyclic bar construction $sk_{n-1} B^{cy} (A;M)$ (which exists even if $A$ is only $A_n$). In particular, a $\mathcal{W}$-trace is corepresented by $B^{cy}(A;M)$.
\end{observation}

If $R(1)=I$, we will call an $R$-trace on $(A,M)$ which restricts to $f : R(1) \otimes M \cong M \sar B$ an $R$-trace extending $f$.

We picture a map $R(S) \otimes M \otimes  A^{\otimes S-\{0\}} \sar B$ as a tree
\begin{equation}
\epsfig{file=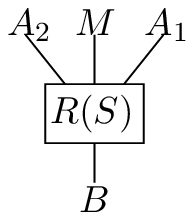}
\end{equation}
where the leaves are cyclically ordered and the leaf labeled $M$ is the basepoint. Here $A_i$ is the copy of $A$ labeled by $i$ in $S=\{0,1,\ldots,n\}$.

\begin{figure} [h]
\epsfig{file=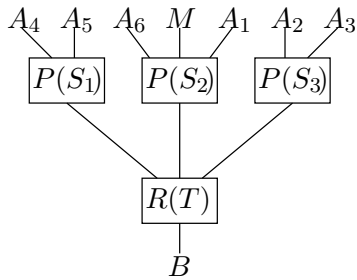}
\caption{A partition of $S$ into $3$ sets.} \label{trace_tree}
\end{figure}

Given a map $f:S \sar T$ in $^0 \! \bD C$, the commutativity of Diagram \ref{trace_comm} says that the two ways to interpret the operation in a diagram such as the one in Figure \ref{trace_tree} agree.

\subsection*{Cotraces}
\begin{figure} [h]
\epsfig{file=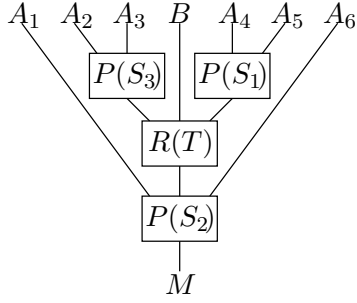}
\caption{Rerooting the tree in Figure \ref{trace_tree}.} \label{cotrace_tree}
\end{figure}

We can also reverse the role of $B$ and $M$, and consider the functor $\tilde{\E}_{B,A,M} : {^0 \! \bD C_P^{op}} \sar \C$ defined by $S \mapsto Hom(B \otimes A^{\otimes S-\{0\}}, M)$. In this case, we interpret the operation corresponding to a tree by first rerooting the tree, making the basepoint leaf the new root. For example, after rerooting the tree in Figure \ref{trace_tree}, it looks like the tree in Figure \ref{cotrace_tree}.
Note that the cyclic ordering in Figure \ref{trace_tree} has been replaced by a linear ordering of the $A$-factors.

\begin{definition}
An $R$-cotrace on $(A,M)$ with source $B$ is a natural transformation $R \sar \tilde{\E}_{B,A,M}$ of functors.
\end{definition}

Giving an $R$-cotrace is equivalent to giving maps $R(S) \otimes B \otimes A^{\otimes S-\{0\}} \sar M$ for each $S \in {^0 \! \bD C}$ such that the diagram
\begin{equation} \label{cotrace_comm}
\xymatrix{
{R}(T) \otimes P[f] \otimes B \otimes  A^{\otimes S-\{0\}} \ar[r] \ar[d] & R(S) \otimes B \otimes A^{\otimes S-\{0\}} \ar[dd] \\
P(f^{-1}(0)) \otimes R(T) \otimes B \otimes A^{\otimes T-\{0\}} \otimes A^{\otimes f^{-1}(0)-\{0\}} \! \! \! \! \! \! \! \! \ar[d] & \\
P(f^{-1}(0)) \otimes M \otimes A^{\otimes f^{-1}(0)-\{0\}} \ar[r] & M}
\end{equation}
commutes for any $f:S \sar T$ in $^0 \! \bD C$. Here the top left vertical map is obtained by writing $A^{\otimes S-\{0\}}$ as $A^{\otimes S-f^{-1}(0)-\{0\}} \otimes A^{\otimes f^{-1}(0)-\{0\}}$ and then using using the maps $P(f^{-1}(t)) \otimes A^{\otimes f^{-1}(t)} \sar A$ for each $t \neq 0$. The fact that this diagram has an extra term corresponds to the fact that rerooting a tree with $2$ levels yields a tree with $3$ levels, as in Diagram  \ref{cotrace_tree}.

A cotrace is represented by a certain object. This is dual to the notion of a trace, but we have to use an extra adjunction, so we present it as a lemma.

\begin{lemma}
$R$-cotraces are represented by $Hom_{{}^0 \! \bD C_P}(R,C^\bullet_{cy}(A;M))$.
\end{lemma}

\begin{proof}
Giving an $R$-cotrace of $B$ into $(A,M)$ is equivalent to giving maps $R(S) \otimes B \otimes A^{\otimes S-\{0\}} \sar M$ which satisfy certain coherence relations. But giving maps $R(S) \otimes B \otimes A^{\otimes S-\{0\}} \sar M$ is equivalent to giving maps $B \sar Hom(R(S) \otimes A^{\otimes S-\{0\}},M)$, and the coherence conditions translate into the conditions for equalizing the maps defining $Hom_{{}^0 \! \bD C_P}(R,$ $C^\bullet_{cy}(A;M))$.
\end{proof}

Again we single out the associahedra and cyclohedra case:

\begin{observation} \label{Wncotraceob}
With $P=\K_n$, $1 \leq n \leq \infty$, a $\mathcal{W}_n$-cotrace is represented by the partial cyclic cobar construction $Tot^{n-1} C_{cy}(A;M)$ (which exists even if $A$ is only $A_n$).
\end{observation}

If we restrict a trace map $W_n \times M \times A^{n-1} \sar B$ extending $f:M \sar B$ to one of the $n$ faces of $W_n$ of the form $K_n$, we get the map 
\begin{equation}
K_n \times M \times A^{n-1} \arrow{1 \times t} K_n \times A^{i-1} \times M \times A^{n-i} \lar M \arrow{f} B,
\end{equation}
where $t$ is the cyclic permutation of $M \times A^{n-1}$ which puts the last $i-1$ factors of $A$ at the beginning and the second map is one of the maps defining the bimodule structure on $M$.

If we restrict a cotrace map $W_n \times B \times A^{n-1} \sar M$ extending $f:B \sar M$ to one of the $K_n$-faces, we get the map
\begin{equation}
K_n \times B \times A^{n-1} \arrow{1 \times \tilde{t}} K_n \times A^{i-1} \times B \times A^{n-i} \arrow{1 \times f \times 1} K_n \times A^{i-1} \times M \times A^{n-i} \lar M,
\end{equation}
where $\tilde{t}$ is the permutation of $B \times A^{n-1}$ placing $B$ in the $i$'th position. Note that in this case there is no cyclic permutation of the factors, the cyclic ordering of $M \times A^{n-1}$ has been replaced by a linear ordering of the $A$-factors in $B \times A^{n-1}$.

                                                                                                       \

\bibliographystyle{plain}
\bibliography{b}
\vspace{12pt}
\noindent
Department of Mathematics, University of Chicago
\newline
5734 S University Ave
\newline
Chicago, IL 60637
\newline
Email: vigleik@math.uchicago.edu

\end{document}